\newcommand\alp{\alpha}

		\newcommand\Del{\Delta}

\newcommand\lam{\lambda}		
\newcommand\sig{\sigma}		

\newcommand\ome{\omega}		\newcommand\Ome{\Omega}

\newcommand\calM{{\cal{M}}}

\newcommand\RR{\Bbb{R}}

\newcommand\ZZ{\Bbb{Z}}

\newcommand\CC{\Bbb{C}}

\newcommand\nek{,\ldots,}
\newcommand\sdp{\times \hskip -0.3em {\raise 0.3ex
\hbox{$\scriptscriptstyle |$}}} 

\newcommand\End{\operatorname{End\,}}

\newcommand\Ker{\operatorname{Ker}}

\newcommand\rk{\operatorname{rk}}

\newcommand\Tr{\operatorname{Tr}}

\newcommand\om{{\overline{m}}}

\newcommand\tilchi{{\widetilde{\chi}}}

\newcommand\ten{\otimes}

\documentstyle[amscd,amssymb,verbatim,12pt]{amsart}

\theoremstyle{plain}
\newtheorem{Thm}[subsection]{Theorem}
\newtheorem{Cor}[subsection]{Corollary}
\newtheorem{Lem}[subsection]{Lemma}
\newtheorem{Prop}[subsection]{Proposition}

\theoremstyle{definition}
\newtheorem{Def}[subsection]{Definition}

\theoremstyle{remark}

\newtheorem{Rem}[subsection]{Remark}

\numberwithin{equation}{section}
\renewcommand{\rm}{\normalshape}

\newif\ifShowLabels
\ShowLabelstrue
\newdimen\theight
\def\TeXref#1{%
	\leavevmode\vadjust{\setbox0=\hbox{{\tt
		\quad\quad  {\small \rm #1}}}%
	\theight=\ht0
	\advance\theight by \lineskip
	\kern -\theight \vbox to
	\theight{\rightline{\rlap{\box0}}%
	\vss}%
	}}%

\ShowLabelsfalse

\renewcommand{\sec}[2]{\section{#2}\label{S:#1}%
	\ifShowLabels \TeXref{{S:#1}} \fi}
\newcommand{\ssec}[2]{\subsection{#2}\label{SS:#1}%
	\ifShowLabels \TeXref{{SS:#1}} \fi}

\newcommand{\refs}[1]{Section ~\ref{S:#1}}
\newcommand{\refss}[1]{Section ~\ref{SS:#1}}
\newcommand{\reft}[1]{Theorem ~\ref{T:#1}}

\newcommand{\refe}[1]{\eqref{E:#1}}

\newenvironment{thm}[1]%
	{ \begin{Thm} \label{T:#1}  \ifShowLabels \TeXref{T:#1} \fi }%
	{ \end{Thm} }

\newcommand{\th}[1]{\begin{thm}{#1} }
\renewcommand{\eth}{\end{thm} }

\newenvironment{lemma}[1]%
	{ \begin{Lem} \label{L:#1}  \ifShowLabels \TeXref{L:#1} \fi }%
	{ \end{Lem} }
\newcommand{\lem}[1]{\begin{lemma}{#1}}
\newcommand{\elem}{\end{lemma}}

\newenvironment{propos}[1]%
	{ \begin{Prop} \label{P:#1}  \ifShowLabels \TeXref{P:#1} \fi }%
	{ \end{Prop} }
\newcommand{\prop}[1]{\begin{propos}{#1}}
\newcommand{\eprop}{\end{propos}}

\newenvironment{corol}[1]%
	{ \begin{Cor} \label{C:#1}  \ifShowLabels \TeXref{C:#1} \fi }%
	{ \end{Cor} }
\newcommand{\cor}[1]{\begin{corol}{#1}}
\newcommand{\ecor}{\end{corol}}

\newenvironment{defeni}[1]%
	{ \begin{Def} \label{D:#1}  \ifShowLabels \TeXref{D:#1} \fi }%
	{ \end{Def} }
\newcommand{\defe}[1]{\begin{defeni}{#1}}
\newcommand{\edefe}{\end{defeni}}

\newenvironment{remark}[1]%
	{ \begin{Rem} \label{R:#1}  \ifShowLabels \TeXref{R:#1} \fi }%
	{ \end{Rem} }
\newcommand{\rem}[1]{\begin{remark}{#1}}
\newcommand{\erem}{\end{remark}}

\newcommand{\eq}[1]%
	{ \ifShowLabels \TeXref{E:#1} \fi
	   \begin{equation} \label{E:#1} }
\newcommand{\eeq}{ \end{equation} }

\newcommand{\prf}{ \begin{pf} }
\newcommand{\epr}{ \end{pf} }

\renewcommand{\b}{\bullet}
\newcommand{\h}[1]{\text{\( H^{#1}(M,F)\)}}
\newcommand{\hb}{\h{\b}}
\newcommand{\dh}[1]{\text{\( \det H^{#1}(M,F)\)}}
\newcommand{\dhb}{\text{\( \det H^{\b}(M,F)\)}}

\newcommand{\ts}{\text{\( C^{\b}(W^u,F)\)}}
\newcommand{\dts}{\text{\( \det C^{\b}(W^u,F)\)}}

\renewcommand{\om}{\text{\( \Ome^{\b}(M,F)\)}}

\newcommand{\omot}{\text{\( \Ome^{\b,[0,1]}_t(M,F)\)}}

\newcommand{\rsm}{\text{\( \|\cdot\|^{RS}_{\det \h{\b}}\)}}
\newcommand{\rsmt}{\text{\( \|\cdot\|^{RS}_{\det \h{\b},t}\)}}

\newcommand{\mm}{\text{\( \|\cdot\|^\calM_{\det \h{\b}}\)}}
\newcommand{\mmt}{\text{\( \|\cdot\|^\calM_{\det \h{\b},t}\)}}

\newcommand{\rsmb}{\text{\( |\cdot|^{RS}_{\det \h{\b}}\)}}
\newcommand{\rsmbt}{\text{\( |\cdot|^{RS}_{\det \h{\b},t}\)}}

\newcommand{\mmb}{\text{\( |\cdot|^\calM_{\det \h{\b}}\)}}

\newcommand{\mt}{\text{\( \rho^\calM\)}}
\newcommand{\rst}{\text{\( \rho^{RS}\)}}
\newcommand{\rstt}{\text{\( \rho^{RS}(t)\)}}

\newcommand{\gtm}{\text{\( g^{TM}\)}}
\newcommand{\gft}{\text{\( g^F_t\)}}
\newcommand{\gf}{\text{\( g^F\)}}

\newcommand{\nf}{\text{\( \nabla f\)}}
\newcommand{\ntm}{\text{\( \nabla^{TM}\)}}
\newcommand{\etm}{\text{\( e(TM,\ntm)\)}}
\newcommand{\mqf}{\text{\( \psi(TM,\ntm)\)}}
\newcommand{\thet}{\text{\( \theta(F,\gf)\)}}

\newcommand{\p}{\text{\( \partial\)}}

\begin{document}


\title{Witten deformation of Ray-Singer analytic torsion}
\author{Maxim Braverman}
\address{School of Mathematical Sciences\\
Tel-Aviv University\\
Ramat-Aviv 69978, Israel}
\email{maxim@@math.tau.ac.il}


\date{August, 1994}

\maketitle

\begin{abstract}
 Let $F$ be a flat vector bundle over a compact Riemannian
manifold $M$ and let
$f:M\to \RR$ be a self-indexing Morse function.
Let $g^F$ be a smooth Euclidean
metric on $F$, let $g^F_t=e^{-2tf}g^F$  and let $\rho^{RS}(t)$
be the Ray-Singer analytic torsion of $F$ associated to the
metric $g^F_t$.
Assuming that $\nabla f$ satisfies the Morse-Smale
transversality conditions,
we provide an asymptotic expansion for $\log \rho^{RS}(t)$ for
$t\to\infty$
of the form $a_0+a_1t+b\log\left(\frac t\pi\right)+o(1)$. We
present explicit
formulae for coefficients $a_0,a_1$ and $b$. In particular, we show that
$b$ is a half integer.
\end{abstract}


\sec{introd} {Introduction}

\ssec{rst}{The Ray-Singer analytic torsion} Let $M$ be
a compact manifold of dimension $n$ and let $F$
be a flat vector bundle on $M$. Let \gf\ and \gtm\ be smooth metrics
on $F$ and $TM$ respectively.

In \cite{rs} Ray and Singer introduced a numerical invariant
of these data which is called the {\em Ray-Singer analytic torsion}
of $F$ and which we shall denote by \rst.

\ssec{wd}{The Witten deformation}Suppose $f:M\to \RR$ is a Morse
function.  For $t>0$, we denote by \gft\  the smooth metric on $F$
\eq{gft}
   	\gft=e^{-2tf}\gf.
\end{equation}
Let $\rstt$ be the Ray-Singer torsion on $F$ associated to the
metrics $\gft$ and $\gtm$.

Denote by \nf\ the gradient vector
field of $f$ with respect to the metric \gtm. Let $B$ be the
finite set of zeroes of \nf.

We shall assume that the following  conditions are satisfied
(cf. \cite[page 5]{bfk3}):
\begin{enumerate}
   	\item $f:M\to \RR$ is a self-indexing Morse function (i.e.
	      $f(x)=index(x)$ for any critical point $x$ of $f$).
	\item The gradient vector field \nf\ satisfies the Smale
	      transversality conditions \cite{sm1,sm2} (for any two
	      critical points $x$ and $y$ of $f$ the stable manifold
	      $W^s(x)$ and the unstable manifold $W^u(y)$,
	      with respect to
	      \nf, intersect transversally).
  	\item For any $x\in B$, the metric \gf\ is flat near $B$ and
  	      there is a system of coordinates
	     $y=(y^1\nek y^n)$ centered at $x$ such that near $x$
	\end{enumerate}
	\eq{met-as}
	  \gtm=\sum_{i=1}^n|dy^i|^2, \quad
	       f(y)=f(x)-\frac12\sum_{i=1}^{index(x)}|y^i|^2+
		\frac12\sum_{i=index(x)+1}^n|y^i|^2.
	\end{equation}

\ssec{as-exp}{Asymptotic expansion of the torsion} Burghelea,
{}Friedlander and Kappeler (\cite{bfk3}) have shown that
the function $\log\rstt$ has
asymptotic expansion for $t\to \infty$ of the form
\eq{bfk}
	\log \rstt=
	\sum_{j=0}^{n+1}a_jt^j+b\log t + o(1).
\end{equation}
The coefficient $a_0$ is calculated in \cite{bfk3}\ in terms
of the parametrix of the Laplace-Beltrami operator.

In the present paper
we shall   calculate all coefficients in the asymptotic
expansion \cite{bfk3}. In fact, we shall show
that the coefficients $a_j=0$ for $j>1$ and
the coefficient $b$ is a half integer.

\ssec{mainres}{}To formulate our result, we need to introduce some
 notation (cf. \cite{bz1}).

Let  \ntm\ be the Levi-Civita connection on $TM$ corresponding to the
metric \gtm, and let \etm\ be the associated representative
of the Euler class of $TM$ in Chern-Weil theory.

Let \mqf\ be the Mathai-Quillen (\cite{mq}) $n-1$ current on $TM$
(see also \cite[Section 3]{bgs}\ and \cite[Section IIId]{bz1}).

Let $\nabla^F$ be the flat connection on $F$ and
let \thet\ be the 1-form on $M$ defined by (cf. \cite[Section IVd]{bz1})
\eq{mqform}
  \thet=\Tr\big[(\gf)^{-1}\nabla^F\gf\big].
\end{equation}
Set
\eq{notation}
  \begin{aligned}
	\chi(F)&=\sum_{i=0}^n(-1)^i\dim H^i(M,F),\\
	\chi'(F)&=\sum_{i=0}^n(-1)^ii\dim H^i(M,F).
  \end{aligned}
\end{equation}

Let \mt\ be the torsion of the Thom-Smale complex (cf. \refs{milmet}).

\th{mainth}
  The function $\log\rstt$ admits an asymptotic expansion for
  $t\to\infty$  of the form
  \eq{mainth}
      \log\rstt=a_0 + a_1t + b\log \left(\frac t\pi\right) + o(1),
  \end{equation}
  where the coefficients $a_0,a_1$ and $b$ are given by the formulas
  \begin{gather} \label{E:a_0}
	a_0=\log\mt -\frac12\int_M\thet(\nf)^*\mqf;\\   \label{E:a1}
      a_1=-\rk(F)\int_Mf\etm+\chi'(F);  \\
\intertext{and}      \label{E:b}
      b=\frac n4 \chi(F)-\frac12\chi'(F).
  \end{gather}
\eth

\rem{halfint} Note that $\chi(F)=0$ if $n$ is odd. Hence, \refe{b}
implies that the coefficient $b$ is a half integer for any $n$.
\erem
\ssec{methpr}{The method of the proof}  Our method is completely
different from that of \cite{bfk3}. In \cite{bfk3}\ the  asymptotic
expansion \refe{bfk}\ is proved by direct analytic arguments and,
then is applied to get a new proof of the
Ray-Singer conjecture \cite{rs}.

In the present paper we use the Bismut-Zhang extension
of this conjecture (\cite{bz1}) in order to obtain
the \reft{mainth}.

\subsection*{Acknowledgments} It is a great pleasure for me
to express  my deep gratitude to Michael Farber for
bringing the papers \cite{bfk2,bfk3} to my attention
and for valuable discussions.

\sec{milmet}{Milnor metric and Milnor torsion}

In this section we follow \cite[Chapter I]{bz1}.
\ssec{detline}{The determinant line of the cohomology}
Let $\h{\b}= \bigoplus_{i=0}^n\h{i}$ be the cohomology
of $M$ with coefficients in $F$ and
let \dh{\b} be the line
\eq{det}
	\dh{\b}=\bigotimes_{i=0}^n \Big(\dh{i}\Big)^{(-1)^i}.
\end{equation}

\ssec{th-sm}{The Thom-Smale complex} Suppose $f:M\to\RR$
is a Morse function satisfying the Smale transversality
conditions \cite{sm1,sm2}\ (for any two
critical points $x$ and $y$ of $f$ the stable manifold
$W^s(x)$ and the unstable manifold $W^u(y)$, with respect to
\nf, intersect transversally).

Let $B$ be the set of critical points of $f$.   If $x\in B$,
we denote by $F_x$  the fiber of $F$ over $x$ and by $[W^u(x)]$
the real line generated by $W^u(x)$. For $0\le i\le n$,
set
\eq{th-sm}
  C^i(W^u,F)=
    \bigoplus_{\begin{Sb} x\in B\\ index(x)=i \end{Sb}}
		[W^u(x)]^*\ten_{\RR} F_x.
\end{equation}
By a basic result of Thom (\cite{thom}) and Smale (\cite{sm2})
(see also \cite[pages 28--30]{bz1}), there is a well defined
linear operators
\[ \p:C^i(W^u,F)\to C^{i+1}(W^u,F), \]
such that  the pair $(C^\b(W^u,F),\p)$ is a complex  and there is
a canonical identification of $\ZZ$-graded vector spaces
$H^\b(C^\b(W^u,F),\p)\simeq H^\b(M,F)$.
By \cite{kmu}\ there is a canonical
isomorphism
\eq{tsc=c}
 	\dh{\b}\simeq \det C^\b(W^u,F).
\end{equation}

\ssec{milmet}{The Milnor metric}
The metric \gf\ on $F$ determines the structure of
Euclidean vector space
on \ts.

\defe{milmet}
The {\em Milnor metric} \mm\ on the line \dh{\b}\ is the metric
corresponding to the obvious metric on \dts\ via the canonical
isomorphism
\refe{tsc=c}.
\edefe

\rem{1}By Milnor \cite[Theorem 9.3]{mi1}, if \gf\ is a flat
metric on $F$ then
the Milnor metric  coincides with the Reidemeister metric defined
through a smooth triangulation of $M$. In this case \mm\
does not depend upon
$F$ and \gtm\ and, hence, is a topological invariant of the flat
Euclidean vector bundle $F$.
\erem

\ssec{miltor}{The Milnor torsion} Let $\p^*$ be the adjoint of $\p$ with
respect to the Euclidean structure on \ts. Using the finite dimensional
Hodge theory, we have the canonical identification
\eq{hodge-ts}
 	H^i(\ts,\p)\simeq \{v\in C^i(W^u,F):\
	\p v=0, \p^* v=0\},\ \ 0\le i\le n.
\end{equation}
As a vector subspace of $C^i(W^u,F)$, the vector space in the right-hand
side of  \refe{hodge-ts} inherits the Euclidean metric. We denote by
\mmb\ the corresponding metric on \dh{\b}.

The metrics \mm\ and \mmb\ do not coincide in general.
We shall describe the
discrepancy.

Set $\Del=\p\p^*+\p^*\p$ and let $P:\ts\to\Ker \Del$ be the orthogonal
projection. Set $\Pi^\perp=1-\Pi$.

Let $N$ and $\tau$ be the operators on \ts\ acting on
$C^i(W^u,F)\ (0\le i\le n)$ by multiplication by $i$
and $(-1)^i$ respectively.
If $A\in \End(\ts)$, we define the supertrace $\Tr_s[A]$ by the formula
\eq{suptr}
	\Tr_s[A]=\Tr[\tau A].
\end{equation}
{}For $s\in\CC$, set
\[
\eta^\calM(s)=-\Tr_s\big[N(\Del)^{-s}\Pi^\perp\big].  \]

\defe{miltor}
The {\em Milnor torsion} is the number
\eq{miltor}
	\mt=\exp\Big(\frac12\frac{d\eta^\calM(0)}{ds}\Big).
\end{equation}
\edefe
The following result is proved in \cite[Proposition 1.5]{bgs}
\eq{mm-mmb}
	\mm=\mmb \cdot\mt.
\end{equation}

\ssec{wit-mil}{Deformation of Milnor metric} The metric \mm\
depends on the metric \gf.  Let $\gft=e^{-2tf}\gf$
and let \mmt\ be the corresponding Milnor metric.
Set
\eq{tilchi}
	\tilchi'(F)= \rk(F)\sum_{x\in B}(-1)^{index(x)}index(x).
\end{equation}
As  $f$ is a self-indexing Morse function
\[
	\rk(F) \sum_{x\in B}(-1)^{index(x)}f(x)=\tilchi'(F).
\]
Obviously,
\eq{wit-mil}
	\mmt= e^{-t\tilchi'(F)}\cdot\mm.
\end{equation}


\sec{rsmet}{Ray-Singer metric and Ray-Singer torsion}

\ssec{badrsmet}{$L_2$ metric on the determinant line}
Let $(\om,d^F)$ be the de Rahm complex of the smooth sections of
$\bigwedge(T^*M)\ten F$ equipped with the coboundary operator $d^F$. The
cohomology of this complex is canonically isomorphic to \hb.

Let $\ast$ be the Hodge operator associated to the metric \gtm. We equip
\om\ with the inner product
\eq{dr-inpr}
	\langle\alp,\alp'\rangle_{\om}=\int_M\langle\alp\land\ast
	\alp'\rangle_{\gf}.
\end{equation}
By Hodge theory, we can identify \hb\ to the corresponding
harmonic forms
in \om. These forms inherit the Euclidean product
\refe{dr-inpr}. Thus the line \dhb\ inherits a
metric \rsmb, which is also called the $L_2$ metric.

\ssec{rstor}{The Ray-Singer torsion}
Let $d^{F*}$ be the formal adjoint of $d^F$ with respect
to the metrics \gtm\
and \gf.

Set $\Del=d^Fd^{F\ast}+d^{F\ast}d^F$ and let
$\Pi:\om\to\Ker \Del$ be the orthogonal
projection. Set $\Pi^\perp=1-\Pi$.

Let $N$ be the operator defining the $\ZZ$-grading of \om, i.e. $N$ acts
on $\Ome^i(M,F)$ by multiplication by $i$.

If an operator $A:\om\to \om$ is trace class,
we define its supertrace $\Tr_s[A]$ as in \refe{suptr}.

{}For $s\in\CC$, set
\[  \eta^{RS}(s)=-\Tr_s\big[N(\Del)^{-s}\Pi^\perp\big].  \]

By a result of Seeley \cite{se}, $\eta^{RS}(s)$ extends to a meromorphic
function of $s\in\CC$, which is holomorphic at $s=0$.

\defe{rstor}
The {\em Ray-Singer torsion \/} is the number
\eq{rstor}
	\rst=\exp\left(\frac12\frac{d\eta^{RS}(0)}{ds}\right).
\end{equation}
\edefe

\ssec{rsmet}{The Ray-Singer metric} We now remind the
following definition
(cf. \cite[Definition 2.2]{bz1}):
\defe{rsmet}
  The {\em Ray-Singer metric \/} \rsm\ on the line \dhb\ is the product
\eq{rsmet}
  \rsm=\rsmb\cdot\rst.
\end{equation}
\edefe

\rem{}When $M$ is odd dimensional, Ray and Singer
\cite[Theorem 2.1]{rs}\
proved that the metric \rsm\ is a topological invariant, i.e. does not
depend on the metrics \gtm\ or \gf.  Bismut and Zhang
\cite[Theorem 0.1]{bz1}\
described explicitly the dependents of \rsm\ on \gtm\ and \gf\ in
the  case when $\dim M$ is odd.
\erem

\ssec{mm=rm}{Bismut-Zhang theorem}
%
%
%
%
Let  \ntm\ be the Levi-Civita connection on $TM$ corresponding to the
metric \gtm, and let \etm\ be the associated representative
of the Euler class of $TM$ in Chern-Weil theory.

Let \mqf\ be the Mathai-Quillen (\cite{mq}) $n-1$ current on $TM$
(see also \cite[Section3]{bgs}\ and \cite[Section IIId]{bz1}).

Let $\nabla^F$ be the flat connection on $F$ and
let \thet\ be the 1-form on $M$ defined by (cf. \cite[Section IVd]{bz1})
\eq{mqform'}
  \thet=\Tr\big[(\gf)^{-1}\nabla^F\gf\big].
\end{equation}

Now we remind the  following theorem by Bismut and Zhang
\cite[Theorem 0.2]{bz1}.
\begin{Thm}[Bismut-Zhang]\label{T:bz}
The following identity holds
   \eq{bz}
	  \log\left(\frac\rsm\mm\right)^2=
		-\int_M\thet(\nf)^*\mqf.
   \end{equation}
\end{Thm}

\ssec{dfrm-bz}{Dependence on the metric}
The metrics \rsm\ and \mm\ depend,
in general, on the metric \gf. Let $\gft=e^{-2tf}\gf$ and let \rsmt\
and \mmt\ be the Ray-Singer and Milnor metrics on \dhb\
associated to the
metrics \gft\ and \gtm.

By \cite[Theorem 6.3]{bz1}
\begin{multline}\label{E:def-lhs}
	\int_M\theta(F,\gft)(\nf)^*\mqf=\\
	     \int_M\thet(\nf)^*\mqf+2t\rk(F)\int_Mf\etm -2t\tilchi'(F).
\end{multline}
{}From \refe{wit-mil}, \refe{bz} and \refe{def-lhs}, we get
\begin{multline}\label{E:wit-milrs}
	\log\left(\frac\rsmt\mm\right)^2=\\
	-\int_M\thet(\nf)^*\mqf- 2t\rk(F)\int_Mf\etm .
\end{multline}

 \sec{mainres}{The main result}

In this section we prove \reft{mainth}, which we restate for
convenience.
\th{mainth'}
  The function $\log\rstt$ admits an asymptotic expansion for
  $t\to\infty$  of the form
  \eq{mainth'}
      \log\rstt=a_0 + a_1t + b\log \left(\frac t\pi\right) + o(1),
  \end{equation}
  where the coefficients $a_0,a_1$ and $b$ are given by the formulas
  \begin{gather} \label{E:a0'}
	a_0=\log\mt -\frac12\int_M\thet(\nf)^*\mqf;\\   \label{E:a1'}
      a_1=-\rk(F)\,\int_Mf\etm+\chi'(F);  \\
\intertext{and}      \label{E:b'}
      b=\frac n4 \chi(F)-\frac12\chi'(F).
  \end{gather}
\eth
%
%
\prf
{}For each $t>0$ we equip \om\ with the inner product
\eq{in-pr-t}
	\langle\alp,\alp'\rangle_{\om,t}=\int_M\langle\alp\land\ast
	\alp'\rangle_{\gft}.
\end{equation}
and we denote by  \rsmbt\  the
$L_2$ metric on \dhb\ (cf. \refss{badrsmet}) associated to this inner
product.

{}From \refe{mm-mmb}, \refe{rsmet} and \refe{wit-milrs}, we get
\begin{multline}\label{E:pr1}
	\log \rstt= -\frac12\int_M\thet(\nf)^*\mqf-\\
	t\rk(F)\,\int_Mf\etm+ \log\mt+\log\left(\frac\mmb\rsmbt\right).
\end{multline}

Let $d^{F\ast}_t$ be the formal adjoint of $d^F$ with respect to the
inner product \refe{in-pr-t}. Set
$\Del_t=d^Fd^{F\ast}_t+d^{F\ast}_td^F$.

Let \omot\  be the direct sum  of the eigenspaces
of $\Del_t$ associated to eigenvalues $\lam\in[0,1]$.
The pair $(\omot,d^F)$ is a subcomplex of $(\om,d^F)$.

We denote by $\|\cdot\|_{\om}$ the norm on \om\ determined by inner product
\refe{in-pr-t}, and by $\|\cdot\|_{\ts}$ the norm on \ts\ determined by \gf\
(cf. \refss{milmet}).

In the sequel, $o(1)$ denotes an element of $\End\big(\ts\big)$ which
preserves the $\ZZ$-grading and is $o(1)$ as $t\to\infty$.

It is shown in \cite[Theorem 6.9]{bz2} that if $t>0$ is large enough,
there exists an isomorphism
\[
	e_t:\ts\to\omot
\]
of $\ZZ$-graded  Euclidean vector spaces such that
\begin{equation}\label{E:et}
	e_t^*e_t=1+o(1).
\end{equation}
By \cite[Theorem 6.11]{bz2}, for any $t>0$ there is
a quasi-isomorphism of complexes
\[ P_t:\Big(\omot,d^F\Big)\to \Big(\ts,\p\Big), \]
which induces the canonical isomorphism
\begin{equation}\label{E:isom}
	 \hb\simeq H^\b(\omot,d^F)\simeq H^\b(\ts,\p)
\end{equation}
and such that
\begin{equation}\label{E:pt-et}
	P_te_t=e^{tN}\left(\frac t\pi \right)^{n/4-N/2}\big(1+o(1)\big).
\end{equation}
Here $e^{tN}\left(\frac{t}\pi\right)^{n/4-N/2}$
denotes the operator on \ts\ acting on $C^i(W^u,F)$ by multiplication
by $e^{ti}\left(\frac{t}\pi\right)^{n/4-i/2}$.

It follows from \refe{pt-et}, that, for $t>0$
large enough, $P_t$ is  one to one.

{}From \refe{et}, \refe{pt-et} we get
\eq{pp*}
	P_tP_t^*=e^{2tN}\left(\frac{t}\pi\right)^{n/2-N}\big(1+o(1)\big).
\end{equation}

{}Fix $\sig\in\h{i}\ (0\le i\le n) $ and let $\ome_t\in \Ker \Del_t$
be the harmonic form representing $\sig$.

Let $\Pi:\ts\to\Ker(\p\p^*+\p^*\p)$ be the orthogonal projection.
Then $\Pi P_t\ome_t\in C^i(W^u,F)$ corresponds to $\sig$ via the canonical
isomorphisms \refe{hodge-ts}, \refe{isom}.

Obviously,
\eq{ker}
	P_t\ome_t\in \Ker\p,\qquad
	e^{2ti}\left(\frac{t}\pi\right)^{n/2-i}\big(P_t^*\big)^{-1}
	\ome_t\in\Ker \p^*.
\end{equation}
By \refe{pp*}, we get
$e^{2ti}\left(\frac{t}\pi\right)^{n/2-i}\big(P_t^*\big)^{-1}\ome_t=
\big(1+o(1)\big)P_t\ome_t$. Then \refe{ker} implies
\eq{pi1}
	\big\|\Pi P_t\ome_t\big\|_{\ts}=
	  \big\|P_t\ome_t\big\|_{\ts}\big(1+o(1)\big).
\end{equation}
{}From \refe{pp*}, \refe{pi1} we obtain
\eq{final}
	\big\|\Pi P_t\ome_t\big\|_{\ts}=
	e^{ti}\left(\frac{t}\pi\right)^{n/4-i/2}
	\big\|\ome_t\big\|_{\om,t}\big(1+o(t)\big).
\end{equation}
It follows from \refe{final} and from the definitions of the metrics
\mmb,  \rsmbt\ that
\eq{pr3}
	\log\left(\frac\mmb\rsmbt\right)=
	t\chi'(F)+\Big(\frac n4\chi(F)-\frac12\chi'(F)\Big)
	\log\left(\frac t\pi \right)+o(1).
\end{equation}
{}From \refe{pr1}, \refe{pr3} we get
\begin{multline}\label{E:pr4}
	\log \rstt=\\
	   -\frac12\int_M\thet(\nf)^*\mqf-t\rk(F)\,\int_Mf\etm+\\
	     \log\mt+t\chi'(F)+
	   \Big(\frac n4\chi(F) - \frac12\chi'(F)\Big)
		\log\left(\frac t\pi \right)+o(1).
\end{multline}
The proof of \reft{mainth'} is completed.
\epr


\end{document}